\documentclass[12pt]{article}\usepackage{amsmath}
\begin{document}
\bibliographystyle{unsrt}\begin{center}{\large\bf Almost Everywhere Flatness
of a 3-Space with a Loop-Based Wormhole}\\[5mm]
Z.~Ya.~Turakulov\\{\it Inter-University Centre for Astronomy and Astrophysics,
Post Bag 4,\\ Ganeshkhind, Pune~411~007, India\\ and Institute of Nuclear
Physics\\ Ulugbek, Tashkent 702132, Uzbekistan\\ (e-mail:turakulov@yahoo.com)}
Comments: LaTeX, 6 pages\\ MSC-class: 53C25 (Primary); 53A99 (Secondary)
\end{center}\begin{abstract}

A particular Riemannian metric which originally has been obtained for a
well-known coordinate system in the Euclidean 3-space, is shown to specify,
in fact, a manifold with boundary. There are two ways to make the manifold
complete. One is to identify two halves of the boundary that turns the
manifold into Euclidean 3-space as it was done originally. Another is to
identify boundaries of two copies of this manifold, that yields a complete
manifold which consists of two copies of Euclidean 3-space connected through a
round disk. In general relativity this kind of connection is called
`loop-based wormhole'. The straightforward calculation of curvature from the
metric specified yields an erroneous result, due which the curvature is zero,
that is impossible because a manifold with this structure cannot be flat. This
paradox is resolved in full correspondence with the generally-accepted
definitions.\end{abstract}\vspace{1cm}\section{Introduction}

Riemannian manifolds are often being specified by explicit form of a metric and
range of each coordinate. Since range of each coordinate is fixed regardless
of the values of all others, they altogether form a $n$-dimensional cuboid in
the space {\bf $R^n$} of $n$-tuples of real numbers. Consequently, this kind of
specification provides a $n$-dimensional cuboid and a metric on it whereas the
standard definition reads that Riemannian manifold is a differentiable manifold
together with a Riemannian metric \cite{k-n,wf}. Evidently, the difference
between a differentiable manifold and its rectangular pattern in {\bf $R^n$}
specified by the coordinate system is that the latter does not possess the same
topology as the earlier. Usually, however, this difference does not reveal
because some details of topology of the desired manifold are known from
general considerations and others follow from the metric itself. It may happen,
however, that some details of structure of a manifold specified this way remain
unknown until a special study is made.

Consider the following example. The metric\begin{equation}ds^2=
(\cosh^2u-\cos^2v)(du^2+dv^2)+\cosh^2u\cos^2vd\varphi^2\end{equation}has been
obtained for oblate spheroidal coordinates in the Euclidean 3-space
\cite{m-p,an}. The ranges of coordinates are given by the inequalities
\begin{equation}0\leq u<\infty,\enskip -\frac{\pi}{2}\leq v\leq\frac{\pi}{2},
\enskip 0\leq\varphi\leq2\pi\end{equation}and the points with $\varphi=0$ and
$\varphi=2\pi$ are identified. Straightforward calculation of the curvature for
this metric yields identically zero result, consequently, from formal point of
view, the metric (1) is flat. Note that this result is obtained disregard of
the range of the coordinate $u$, therefore, it seems that the space with the
metric (1) and ranges of the coordinates given by the inequalities
\begin{equation}-\infty\leq u<\infty,\enskip -\frac{\pi}{2}\leq v\leq
\frac{\pi}{2},\enskip 0\leq\varphi\leq2\pi\end{equation}is flat too. On the
other hand, all non-compact Riemannian manifolds are known and classified in
standard texts \cite{k-n,wf}, but this space is not mentioned in the
classification. Indeed, this manifold apparently differs from all known flat
Riemannian spaces because, at least, it contains two copies of Euclidean
3-space connected through a round disk $u=0$. In this work we show that in
spite of zero result of straightforward calculation of its curvature, this
manifold is non-flat.\section{Euclidean 3-space from the metric (1)}

As was pointed out above, the metric (1) has been obtained for a particular
coordinare system in the Euclidean 3-space. In this section we consider the
inverse problem and will show how to obtain this space from the metric. There
exists a simple transformation from the coordinates $\{u,v,\varphi\}$ to
Cartesian coordinates and, thereby, to Euclidean geometry, therefore, it
suffices to find, where the lower limit of the coordinate $u$ comes from. To
see it, cosider the coordinate surfaces $u=const$.

These surfaces are confocal oblate spheroids. If the lower limit of the
coordinate $u$ is some positive number $u_0$, then the metric (1) with $v$ and
$\varphi$ coordinates ranging as in the inequalities (2) and the inequalities
$0\leq u_0\leq u<\infty$ specify the exterior of the spheroid $u=u_0$ in the
Euclidean 3-space. Note that this manifold has boundary and consider its limit
when $u_0$ tends to zero. In this limit the manifold in question becomes the
exterior of the infinitesimally thin disk of unit radius specified by the
equation $u=0$, however, the boundary remains.

Now, since the two sides of the disk are just circles of unit radius, they can
be identified. As it is done the manifold becomes a Euclidean 3-space, but the
coordinate $v$ is not longer a continuous function on it because of the sign
changing discontinuity created by identificationof the surfaces $u=0,\ v<0$
and $u=0,\ v>0$. Thus, Euclidean 3-space cannot be obtained from the metric (1)
alone with the inequalities (2); its topology ust also be spacified by
identifying the surfaces just mentioned. Otherwise, putting the lower limit of
the coordinate $u$ creates a boundary of the manifold. Consequently, in order to
obtain Euclidean 3-space from the metric (1) we must first create a boundary
cutting the range of the coordinate $u$, then remove it identifying its halves,
therefore a question arises, what happens if these two operations have not been
made. In other words, what manifold is specified by the metric (1) and
inequialities (3)? Evidently, this manifold is complete and, as was pointed out
above, straightforward calculation of curvature made in the coordinates
$\{u,v,\varphi\}$ yields identically zero result. Besides, since the coordinate
$u$ has no limits the metric (1) and inequalities (2) specify the manifold
exhaustively leaving no possibility prescribe any details of its structure.
\section{The meridian surface}

The manifold under consideration has one Killing vector $\partial_\varphi$, and
it is natural to first to explore geometry of a typical surface orthogonal to
it. The coordinate surfaces $\varphi=const$ are incomplete, because under
$u\neq0$ they are half-planes with axis of symmetry as the boundary line. Two
half-planes $\varphi=\varphi_0$ and $\varphi=\varphi_0\pm\pi\leq2\pi$ together
constitute an entire plane orthogonal everywhere to the Killing vector.
Hereafter surfaces defined this way are called the meridian surfaces.

Now, let us return to the exterior of the disk $u=0$, which is a Riemannian
manifold with boundary specified by the metric (1) and inequalities (2)
without identifying the opposite sides of the disk. Its meridian surface is a
plane with straightline cut whose edges are halves of the coordinate line $u=0$,
on which the coordinate $v$ has opposite sign. Note that a meridian surface of
the Riemannian manifold specified by the metric (1) and inequalities (3)
contains two copies of the the surface just considered, on which the coordinate
$u$ has opposite sign, and their boundaries are eliminated by identifying the
edges of the cuts. In other words, the arc $u=0,\ v<0$ of one copy is
identified with the arc $u=0,\ v<0$ of another, and so for $v>0$. Note that on
this surface the coordinate $v$ is a continuous function.

Surprizingly, the manifold defined such an artificial way is, in fact, nothing
but a kind of quadric surface. Indeed, consider an elliptical hyperboloid of
one sheet. In standard Cartesian coordinates it can be expressed by the
equation$$\frac{x^2}{a^2+\lambda}+\frac{y^2}{b^2+\lambda}+\frac{z^2}{\lambda}=
1,\enskip a^2>b^2,\enskip -b^2<\lambda<0.$$Metric of this surace has the form
\begin{equation}ds^2=\frac{(u-\lambda)(u-w)}{4u(a^2+u)(b^2+u)}du^2+
\frac{(w-\lambda)(w-u)}{4w(a^2+w)(b^2+w)}dw^2,\end{equation}where $u>0,\
-a^2<w<-b^2$ \cite{sk,wn}. Note that the section of the hyperboloid with the
plane $z=0$ is the ellipse$$\frac{x^2}{a^2+\lambda}+\frac{y^2}{b^2+\lambda}=1.
$$Now, pass to the limit $b\rightarrow0$, under which the parameter $\lambda$
also takes the zero value. The ellipse in the $z=0$ plane turna into a
straightline segment, consequently, the resulting surface becomes similar to
the meridian surface considered above. The metric (4) becomes$$ds^2=
\frac{u-w}{4u(a^2+u)}du^2-\frac{u-w}{4w(a^2+w)}dw^2,\enskip -a^2<w<-0.$$Now,
the substitutions $u=a^2\sinh^2\zeta,\ w=-a^2\sin^2\eta$ transforms it into the
following:$$ds^2=a^2(\sinh^2\zeta+\sin^2\eta)(d\zeta^2+d\eta^2)$$that under
$a=1$ coincides with the first term in the right-hand side of the equation (1).
Consequently, the meridian surface obtained above by the cut-and-paste
procedure applied to two planes is the limiting case of elliptical hyperboloid
just considered. This fact seems to be surprizing because, on one hand, the
metric obtained for the hyperboloid of this kind is known to be flat whereas
hyperboloids are non-flat manifolds. Curvature of one-sheet hyperboloid is
concentrated mainly in the throat and in the limiting case under consideration
the throat still presents, though in a somewhat unusual form. To explain this
phenomenon, consider the endpoints of the cuts made on both planes alone the
coordinate line $u=0$. It is seen that after identifying the edges of the cuts
these two points become special ones. If we try to introduce polar coordinates
$\{\rho,\theta\}$ on the surface with the pole in  one of them we have to let
the angle $\theta$ range from $0$ to $4\pi$ instead of $2\pi$. Evidently, the
pole is a saddle reduced to a single point, and, since the angular defect is
exactly $-2\pi$, the curvature of the surface in its neighborhood is
$-\delta(\rho)/2\rho$ because measure of surface integration here is
$\rho d\rho d\theta$ and $\theta$ runs from $0$ to $4\pi$. In the coordinates
$\{u,v\}$ this gives:\begin{equation}
R_u{}^v{}_{uv}=-\frac{\delta(u)\delta(v)}{\cosh^2u-\cos^2v}\end{equation}
(note that the coordinate v runs from $-\pi/2$ to $\pi/2$ twice).
\section{A space with a loop-based wormhole}

Now we return to the Riemannian manifold specified by the metric (1) and the
ranges of the coordinates (3) and summarize its geometric properties. The
manifold has rotational symmetry and admits a foliation whose leaves are
orthogonal everywhere to the Killing vector. A typical leaf of the folation is
a limiting case of elliptic hyperboloid of one sheet, whose throat is reduced to
a straightline segment. The manifold consists of two copies of Euclidean 3-space
connected through the round disk formed by rotation of the segment. Any
straightline which crosses the disk passes from one copy of Euclidean 3-space to
another, and all the rest ones lie wholly in only one of them. This phenomenon
and the curvature singularity are known in general relativity as a loop-based
wormhole. It was described in the
book \cite{vr} as follows.

{\it Now cut the two flat three space open aloong the surfaces $S_1$, $S_2$.
Identify $S_1{}^+$ with $S_2{}^-$, and identify $S_1{}^-1$ with $S_2{}^+$.
These are smooth identifications, and in fact there is no discontinuity in the
second fundamental form. Consequently the Riemann tensor is everywhere zero
except possibly at the loops $L_1$, $L_2$ themselves. This is now our model
for a ``loop-based'' wormhole connecting two flat spaces... 

The net result is, that for loop-based womholes constructed using cut and
paste Minkowski spacetimes, the Einstein tensor is
$$G_{\mu\nu}=-\delta^2(x,\Sigma)h_{\mu\nu}$$...}\\
Here the term ``loop'' stands for the edge of the disk and has
coordinates $u=v=0$. Note that the manifold has two axes of symmetry whose
points remain immobile under the actions of the Killing vector and each of them
crosses the disk, thus belongs to both copies of the Euclidean 3-space.
\section{Conclusion}

As the preceding sections made clear, the manifold specified by the metric (1)
and inequalities (3) is flat everywhere but the focal circle of the coordinate
system$\{u,v,\varphi\}$. Though straightforward calculation of the curvature
from the metric yields zero result, the manifold is neither flat nor smooth.
Its structure can be seen from that of a typical surface which is orthogonal
everywhere to the only Killing vector $\partial_\varphi$. This surface is a
limiting case of elliptic hyperboloid of one sheet, whose throat is reduced to
a straightline segment.

Alternatively, this surface can be represented as the result of cut and paste
procedure applied to a pair of planes as described in the section 3. The
endpoints of the segment are saddles in which the curvature has a
$\delta$-function singularity. The structure of the entire manifold, in which
the segment draws a round disk, inherites features of this surface. It consists
of two copies of Euclidean 3-space connected through the (doubled) disk the
same way as two copies of plane are connected via the (doubled) segment, and
the edge of the disk (not doubled) is the support of the curvature of the
manifold. The Riemannian curvature has only one non-zero component given by the
equation (5), which, however, does not follow straightforwardly from the form
of mentric (1).

The form of metric does not expose any non-zero curvature because, on one hand,
the curvature is presented by a $\delta$-function singularity, and, on the other
hand, its support coincides with the singularity of the coordinate system which
is the focal circle $u=v=0$ (note that this curve is the zero of the metric
determinant). It must be emphasized that this example does not contradict the
generally-accepted definition of flat Riemannian manifold \cite{k-n}. Indeed,
the definition requires that the curvature vanishes identically whereas in this
case it does only in a particular coordinate system.
{\bf Acknowledgment:} The author thanks the Third World Academy of Sciences for
financial support and IUCAA for warm hospitality which made this work possible.
\end{document}